\documentclass[11pt]{article}
\usepackage{amsthm, amsmath, amssymb, amsfonts, url, booktabs, tikz, setspace, fancyhdr, bm, mathrsfs}
\usepackage{hyperref}
\usepackage{geometry}
\usepackage{hyperref, enumerate}
\usepackage[shortlabels]{enumitem}
\usepackage[babel]{microtype}
\usepackage[english]{babel}
\usepackage{aliascnt}
\usepackage[capitalise]{cleveref}
\usepackage{comment}
\usepackage{bbm}
\usepackage{csquotes}
\usepackage{mathabx}
\usepackage{tikz}
\usepackage{graphicx}
\usepackage{float}
\usepackage{amsmath}

\counterwithin{figure}{section}

\newtheorem{theorem}{Theorem}[section]
\newaliascnt{prop}{theorem}
\newtheorem{prop}[prop]{Proposition}
\aliascntresetthe{prop}
\newaliascnt{lemma}{theorem}
\newtheorem{lemma}[lemma]{Lemma}
\aliascntresetthe{lemma}
\newaliascnt{cor}{theorem}

\aliascntresetthe{cor}
\newaliascnt{conj}{theorem}
\newtheorem{conj}[conj]{Conjecture}
\aliascntresetthe{conj}
\newaliascnt{claim}{theorem}
\newtheorem{claim}[claim]{Claim}
\aliascntresetthe{claim}

\newaliascnt{fact}{theorem}

\aliascntresetthe{fact}
\crefname{claim}{Claim}{Claims}
\Crefname{claim}{Claim}{Claims}
\crefname{lemma}{Lemma}{Lemmas}
\Crefname{lemma}{Lemma}{Lemmas}
\crefname{cor}{Corollary}{Corollaries}
\Crefname{cor}{Corollary}{Corollaries}
\crefname{prop}{Proposition}{Propositions}
\Crefname{prop}{Proposition}{Propositions}

\theoremstyle{definition}

\newtheorem*{defn-non}{Definition}

\usepackage[linesnumbered, ruled]{algorithm2e}
\SetKwRepeat{Do}{do}{while}
\newtheorem{rmk}[theorem]{Remark}

\newenvironment{poc}{\begin{proof}[Proof of claim]}{\end{proof}}

\usepackage{todonotes}

\newcommand{\cP}{\mathcal{P}}

\newcommand{\cB}{\mathcal{B}}

\newcommand{\cF}{\mathcal{F}}

\newcommand{\cH}{\mathcal{H}}

\newcommand{\cA}{\mathcal{A}}
\newcommand{\cC}{\mathcal{C}}
\newcommand{\cR}{\mathcal{R}}
\DeclareMathOperator{\Sh}{Sh}
\DeclareMathOperator{\VC}{VCdim}
\newcommand{\E}{\mathbb{E}}
\newcommand{\ind}{\mathbbm{1}}
\title{Largest Sperner families with bounded VC-dimension}
\author{Gennian Ge\thanks{School of Mathematical Sciences, Capital Normal University, Beijing, China. Emails: gnge@zju.edu.cn, 2230502181@cnu.edu.cn, 2240501023@cnu.edu.cn. Gennian Ge is supported by the National Key Research and Development Program of China under Grant 2025YFC3409900, the National Natural Science Foundation of China under Grant 12231014, and Beijing Scholars Program.}
\and Hanyu Wang\footnotemark[1]
\and Zixiang Xu\thanks{School of Mathematical Sciences, Zhejiang University, Hangzhou, China. Email: zixiangxu@zju.edu.cn.}
\and Yixuan Zhang\footnotemark[1]
}
\date{}

\begin{document}
\maketitle
\vspace{-1.7em}
\begin{abstract}
A conjecture of Frankl proposed in the 1980s states that for $n\ge2d$, every Sperner family $\cF\subseteq 2^{[n]}$ of VC-dimension at most $d$ has at most $\binom{n}{d}$ members. We prove this conjecture for every fixed $d$ and large $n$, and show that equality holds only when $\cF$ is \(\binom{[n]}{d}\) or \(\binom{[n]}{n-d}.\)
\end{abstract}

\section{Introduction}\label{sec:introduction}

For a finite set $V$, write $\binom{V}{k}=\{A\subseteq V:|A|=k\}$ and $\binom{V}{\le k}=\bigcup_{j=0}^{k}\binom{V}{j}$, and set $[n]=\{1,\ldots,n\}$. We call $\binom{V}{k}$ the \emph{complete layer} of sets of size $k$. A family $\cF\subseteq 2^{[n]}$ is a \emph{Sperner family} if no member contains another. Determining the largest possible size of such a family is one of the earliest problems in extremal set theory. Sperner's theorem~\cite{Sperner1928} states that $|\cF|\le\binom{n}{\lfloor n/2\rfloor}$, which is sharp. 

For $S\subseteq[n]$, the \emph{trace} of $\cF$ on $S$ is $\cF|_{S}=\{A\cap S:A\in\cF\}$. We say that $\cF$ \emph{shatters} $S$ if $\cF|_{S}=2^{S}$, and write $\Sh(\cF)$ for the family of shattered sets. The \emph{Vapnik--Chervonenkis dimension} (VC-dimension for short), denoted by $\VC(\cF)$, is the largest size of a shattered set. The Sauer--Shelah lemma~\cite{Sauer1972,Shelah1972,VapnikChervonenkis1971} states that every family with $\VC(\cF)\le d$ satisfies $|\cF|\le\sum_{i=0}^{d}\binom{n}{i}$. This bound can be attained by $\binom{[n]}{\le d}$. VC-dimension has been studied in several problems in extremal combinatorics. For uniform families, Frankl and Pach~\cite{FranklPach1984} proved that a family of sets of size $d+1$ with VC-dimension at most $d$ has size at most $\binom{n}{d}$ when $n\ge2d$. Subsequent work has given improved upper bounds and new constructions~\cite{AhlswedeKhachatrian1997,ChaoXuYipZhang2025,GeXuYipZhangZhao2026,MubayiZhao2007,TranXu2026,WangXuZhang2025,YangYu2025,ZhaoGe2026b}. Bounded VC-dimension also plays a role in sunflower problems~\cite{BaloghEtAl2025,FoxPachSuk2023,GeWangXuZhang2026} and in graph and Ramsey problems~\cite{FoxPachSuk2019,FoxPachSuk2021,JCTBVC2026,2024GeoToGraph,NguyenScottSeymour2023,2026EHBdforGraph}.

Notice that for $n\ge2d$, the complete layers $\binom{[n]}{d}$ and $\binom{[n]}{n-d}$ are Sperner families of VC-dimension $d$, each of size $\binom{n}{d}$. In 1989, Frankl~\cite{Frankl1989} conjectured that these layers have the largest possible size.

\begin{conj}\label{conj:frankl}
Let $d\ge1$ and $n\ge2d$ be integers, and let $\cF\subseteq2^{[n]}$ be a Sperner family with $\VC(\cF)\le d$. Then $|\cF|\le\binom{n}{d}$.
\end{conj}

This conjecture brings together Sperner's theorem and the Sauer--Shelah lemma. When $n=2d$ or $n=2d+1$, the proposed bound is exactly the bound in Sperner's theorem. Frankl~\cite{Frankl1989} proved the cases $d=1,2$, and Anstee and Sali~\cite{AnsteeSali1997} proved the case $d=3$. Despite subsequent progress on related problems, the cases of Conjecture~\ref{conj:frankl} known in the full range $n\ge2d$ have remained $d\le3$ since 1997. The conjecture was also restated in recent work of Frankl and Wang~\cite[Conjecture~4.8]{2024FourVertexTraces}, Yang and Yu~\cite[Conjecture~4.4]{YangYu2025} and Wang, Xu and Zhang~\cite[Conjecture~3.1]{WangXuZhang2025}.

We prove the conjecture for every fixed $d\ge3$ and large $n$, and determine the extremal families.

\begin{theorem}\label{thm:main}
For every integer $d\ge3$, there exists some \(n_{0}=n_{0}(d)\) such that for every \(n\ge n_{0}\), every Sperner family $\cF\subseteq2^{[n]}$ with $\VC(\cF)\le d$ satisfies
\(
|\cF|\le\binom{n}{d}.
\)
Equality holds if and only if $\cF=\binom{[n]}{d}$ or $\cF=\binom{[n]}{n-d}$.
\end{theorem}

\begin{rmk}\label{rem:threshold}
Our proof shows that it suffices to take $n_{0}(d)=d^{10}\cdot 2^{5d-10}$. The main remaining problem is to establish the conjecture in the full range $n\geq2d$.
\end{rmk}
Frankl~\cite{Frankl1989} provided the general upper bound \(\binom{n}{d}+\binom{n}{d-3}+\binom{n}{d-6}+\cdots\) by a recursive
argument and therefore proved the cases $d\le 2$. Anstee and Sali~\cite{AnsteeSali1997} settled $d=3$
using induction and a characterization of the extremal families
for $d=1,2$. Our proof mainly combines averaging over restrictions,
a counting argument based on the ballot problem in the field of enumerative combinatorics~\cite{RomanSteven2015}, and the
Sauer-Shelah lemma~\cite{Sauer1972,Shelah1972,VapnikChervonenkis1971}.
\section{Proof of the main theorem}\label{sec:proof}

\subsection{Some useful tools}\label{subsec:tools}

We write $A\triangle B=(A\setminus B)\cup(B\setminus A)$ for the symmetric difference of two sets. We first recall two fundamental results related to the main topic of this paper.
\begin{lemma}[\cite{Sperner1928}]\label{lem:Sperner}
Every Sperner family $\cA\subseteq 2^{[n]}$ satisfies $|\cA|\le\binom{n}{\lfloor n/2\rfloor}$.
\end{lemma}

\begin{lemma}[\cite{Sauer1972,Shelah1972,VapnikChervonenkis1971}]\label{lem:sauer}
Let \(d\ge 1\) be an integer and $\cH\subseteq 2^{[n]}$ be a set system with VC-dimension at most $d$. Then $|\cH|\le\sum_{j=0}^{d}\binom{n}{j}$.
\end{lemma}

We also use the result on the ballot problem in enumerative combinatorics~\cite{RomanSteven2015}. 
\begin{lemma}[\cite{RomanSteven2015}]\label{lem:Roman Catalan}
Consider an election with two candidates $A$ and $B$ where $A$ receives $a$ votes and $B$ receives $b$ votes, with $a\ge b$. Assume that all voting histories are equally likely. The probability that candidate A is never behind candidate B is $\frac{a+1-b}{a+1}$.
\end{lemma}

\subsection{Transformation of the problem and the proof of Theorem~\ref{thm:main}}\label{subsec:main}

Let $\cF\subseteq 2^{[n]}$ be a Sperner family with $\VC(\cF)\le d$. For $J\subseteq[n]$ and $U\subseteq[n]\setminus J$, define the \emph{restricted family} $\cF_{J,U}=\{A\subseteq J:A\cup U\in\cF\}$. We use this notation for any family on a specified ground set.
Every restricted family $\cF_{J,U}$ is Sperner. For fixed $J$, each $F\in\cF$ has a unique decomposition
$F=A\cup U$ with $A\subseteq J$ and
$U\subseteq[n]\setminus J$. Therefore, $\sum_{U\subseteq[n]\setminus J}|\cF_{J,U}|=|\cF|$.  

Further, for any $S\subseteq J\subseteq[n]$, define the family
$\cR_{J,S}
=\{U\subseteq[n]\setminus J:S\in\Sh(\cF_{J,U})\}$.
We show that any set shattered by $\cR_{J,S}$ can be joined to $S$ to give
a set shattered by $\cF$.

\begin{lemma}\label{lem:residual}
If $T\in\Sh(\cR_{J,S})$, then $S\cup T\in\Sh(\cF)$.
In particular, $\VC(\cR_{J,S})\le d-|S|$ when $|S|\le d$, and $\cR_{J,S}=\emptyset$ when $|S|>d$.
\end{lemma}

\begin{proof}[Proof of Lemma~\ref{lem:residual}]
Note $T\in\Sh(\cR_{J,S})$ if and only if for any $X\subseteq T$, there exists $U\in\cR_{J,S}$ such that $U\cap T=X$. By the definition of $\cR_{J,S}$, $U\in\cR_{J,S}$ if and only if  $S\in\Sh(\cF_{J,U})$. Then for any $Y\subseteq S$, there exists $A\in\cF$, such that $(A\setminus U)\in\cF_{J,U}$ and $(A\setminus U)\cap S=Y$. So $A\cap (S\cup T)=X\cup Y$, which means $\cF$ shatters $S\cup T$. Note $S$ is disjoint with $T$, so $\VC(\cR_{J,S})\le d-|S|$ when $|S|\le d$. If $\cR_{J,S}$ is nonempty and $|S|>d$, then the empty set can be shattered by $\cR_{J,S}$ and hence $S$ can be shattered by $\cF$. This contradicts $\VC(\cF)\le d$.
\end{proof}

By Lemma~\ref{lem:residual} and Lemma~\ref{lem:sauer}, we have $|\cR_{J,S}|\le\sum_{j=0}^{d-|S|}\binom{n-|J|}{j}$. It is not hard to show that for any family $\cA$ on the finite ground set, $|\cA|\le|\Sh(\cA)|$. Indeed, we can prove it by using the downshift argument in Lemma~\ref{lem:downshift}. Also note $|\cF|
=\sum_{U\subseteq[n]\setminus J}|\cF_{J,U}|$. So we have
\begin{align*}
|\cF|
&\le\sum_{U\subseteq[n]\setminus J}|\Sh(\cF_{J,U})|\\
&=\sum_{S\subseteq J}|\cR_{J,S}|\\
&\le\sum_{s=0}^{\min(|J|,d)}\binom{|J|}{s}\sum_{j=0}^{d-s}\binom{n-|J|}{j}.
\end{align*}
However, $\sum_{s=0}^{\min(|J|,d)}\binom{|J|}{s}\sum_{j=0}^{d-s}\binom{n-|J|}{j}
=\sum_{k=0}^{d}\sum_{s=0}^{\min(|J|,k)}\binom{|J|}{s}\binom{n-|J|}{k-s}
=\sum_{k=0}^{d}\binom{n}{k}$. Thus this gives the original Sauer--Shelah bound.
Splitting the coordinates and counting all shattered sets in this way does not yet improve the bound.

To get the desired bound $\binom{n}{d}$, we assign the same nonnegative weight $w_s$ to every set of size $s$ in $J$. In other words, we want to have the identity $\sum_{s=0}^{\min(|J|,d)}\binom{|J|}{s}w_s
\sum_{j=0}^{d-s}\binom{n-|J|}{j}
=\binom{n}{d}$. That is, $\sum_{j=0}^{d}\binom{n-|J|}{j}
\sum_{s=0}^{\min(|J|,d-j)}\binom{|J|}{s}w_s
=\sum_{j=0}^{d}\binom{n-|J|}{j}\binom{|J|}{d-j}$. So we can choose the weights so that $\sum_{s=0}^{\min(|J|,k)}\binom{|J|}{s}w_s
=\binom{|J|}{k}$ for every $0\le k\le d$. 
With such weights, we have $|J|\ge2d-1$ since $w_s$ is nonnegative for every $0\le s\le d$. 

Put $b_{p,k}=\binom{p}{k}-\binom{p}{k-1}$ for integer $p>0$ and integer $0\le k\le \lfloor p/2\rfloor$. If we take $|J|=2d-1$, then we have the following identity.

\begin{lemma}\label{lem:terminal}
For every integer $n\ge 2d-1$, we have
\[
\sum_{s=0}^{d-1}b_{2d-1,s}\sum_{j=0}^{d-s}\binom{n-2d+1}{j}=\binom{n}{d}.
\]
\end{lemma}

\begin{proof}[Proof of Lemma~\ref{lem:terminal}]
Note we have $\sum_{s=0}^{d-1}b_{2d-1,s}\sum_{j=0}^{d-s}\binom{n-2d+1}{j}
=\sum_{j=0}^{d}\binom{n-2d+1}{j}\sum_{s=0}^{\min(d-j,d-1)}b_{2d-1,s}$ and $\sum_{s=0}^{a}b_{2d-1,s}=\binom{2d-1}{a}$ for any $0\le a\le d-1$. Since $\binom{2d-1}{d-1}=\binom{2d-1}{d-0}$, we have $\sum_{s=0}^{\min(d-j,d-1)}b_{2d-1,s}=\binom{2d-1}{d-j}$ for any $0\le j\le d$. So we have 
\[
\sum_{s=0}^{d-1}b_{2d-1,s}\sum_{j=0}^{d-s}\binom{n-2d+1}{j}
=\sum_{j=0}^{d}\binom{n-2d+1}{j}\binom{2d-1}{d-j}
=\binom{n}{d}.
\]

\end{proof}

An important observation is that $b_{p,k}$ is related to the ballot problem.
For any $J\subseteq[n]$ of odd size $p\ge 3$, let $\boldsymbol{\pi}=(x_{1},\ldots,x_{p})$ be an ordering of $J$. The set $\{x_{s_{1}},\ldots,x_{s_{k}}\}$ with $1\le s_{1}<\cdots<s_{k}\le p$ is called a \emph{ballot set} if $s_{j}\ge 2j$ for every $j\in[k]$. We denote $\cB_{p}(\boldsymbol{\pi})$ for the family of ballot sets, including the empty set. Every member of $\cB_{p}(\boldsymbol{\pi})$ has size at most $\frac{p-1}{2}$.

\begin{lemma}\label{lem:ballot}
For $0\le k\le\frac{p-1}{2}$, there are exactly $b_{p,k}$ members of $\cB_{p}(\boldsymbol{\pi})$ having size $k$.
If $S\in\binom{J}{k}$ and $\boldsymbol{\pi}$ is chosen uniformly at random, then $\mathbb{P}(S\in\cB_{p}(\boldsymbol{\pi}))=b_{p,k}/\binom{p}{k}$.
If $|J|=2d-1$, then for any ordering $\boldsymbol{\pi}$ of $J$, we have 
\[
\sum_{S\in\cB_{2d-1}(\boldsymbol{\pi})}|\cR_{J,S}|\le\binom{n}{d}.
\]
\end{lemma}

\begin{proof}[Proof of Lemma~\ref{lem:ballot}]
 We can see $S$ as a subset of $[p]$ of size $k$ if we see $J$ as $[p]$. For a fixed subset $S$ of size $k$, we call an element in $[p]$ selected if $S$ contains this element, and unselected otherwise. Then by the definition of $\cB_p(\boldsymbol{\pi})$, 
 a nonempty subset of $J$ is contained in $\cB_p(\boldsymbol{\pi})$
if and only if the unselected element is never behind the selected element. By Lemma~\ref{lem:Roman Catalan}, we have $|\cB_p(\boldsymbol{\pi})\cap\binom{J}{k}|=\frac{p-2k+1}{p-k+1}\cdot\binom{p}{k}=b_{p,k}$ and $\mathbb{P}(S\in\cB_p(\boldsymbol{\pi}))=b_{p,k}/\binom{p}{k}$ when the orderings are uniformly distributed.

Now suppose $|J|=2d-1$. Then for any ordering $\boldsymbol{\pi}$ of $J$, by Lemma~\ref{lem:sauer}, we have $\sum_{S\in\cB_{2d-1}(\boldsymbol{\pi})}|\cR_{J,S}|
\le
\sum_{S\in\cB_{2d-1}(\boldsymbol{\pi})}
\sum_{j=0}^{d-|S|}\binom{n-2d+1}{j}$. Since $\left|\cB_{2d-1}(\boldsymbol{\pi})\cap\binom{J}{k}\right|=b_{2d-1,k}$, by Lemma~\ref{lem:terminal}, we have 
\[
\sum_{S\in\cB_{2d-1}(\boldsymbol{\pi})}|\cR_{J,S}|
\le
\sum_{s=0}^{d-1}b_{2d-1,s}
\sum_{j=0}^{d-s}\binom{n-2d+1}{j}
=\binom{n}{d}.
\]

\end{proof}

For any $J\subseteq[n]$ with $|J|=2d-1$ and any ordering $\boldsymbol{\pi}$ of $J$, by the definition of $\cR_{J,S}$, Lemma~\ref{lem:ballot} and double counting the pairs $(U,S)$ with $U\subseteq[n]\setminus J$ and $S\in \Sh(\cF_{J,U})\cap\cB_{2d-1}(\boldsymbol{\pi}) $, we have 
\[
\sum_{S\in\cB_{2d-1}(\boldsymbol{\pi})}|\cR_{J,S}|=
\sum_{U\subseteq[n]\setminus J}
|\Sh(\cF_{J,U})\cap\cB_{2d-1}(\boldsymbol{\pi})|
\le\binom{n}{d}.
\]
Note we want to prove  $|\cF|\le\binom{n}{d}$. Using $|\cF|=\sum_{U\subseteq[n]\setminus J}|\cF_{J,U}|$, we have
\begin{equation}\label{eq:fiber}
\sum_{U\subseteq[n]\setminus J}
\left(
|\Sh(\cF_{J,U})\cap\cB_{2d-1}(\boldsymbol{\pi})|
-|\cF_{J,U}|
\right)
\le\binom{n}{d}-|\cF|.
\end{equation}
For any $\cA\subseteq 2^{J}$ with $J\subseteq[n]$ having odd size $p\ge 3$ and any ordering $\boldsymbol{\pi}$ of $J$, we define $\delta_{p,\boldsymbol{\pi}}(\cA)=|\Sh(\cA)\cap\cB_{p}(\boldsymbol{\pi})|-|\cA|$ and $\overline{\delta}_{p}(\cA)=\E_{\boldsymbol{\pi}}\delta_{p,\boldsymbol{\pi}}(\cA)$, where all $p!$ orderings of $J$ are chosen uniformly. Lemma~\ref{lem:ballot} gives
\begin{equation}\label{eq:defect}
\overline{\delta}_{p}(\cA)=\sum_{k=0}^{(p-1)/2}\frac{b_{p,k}}{\binom{p}{k}}\bigl|\Sh(\cA)\cap\binom{J}{k}\bigr|-|\cA|.
\end{equation}
We can also average over all $J\in\binom{[n]}{p}$. Therefore for odd $3\le p\le n$, we define
\begin{equation}\label{eq:alpha}
    \alpha_p(\cF)
=\frac{1}{\binom{n}{p}}
  \sum_{J\in\binom{[n]}{p}}
  \sum_{U\subseteq[n]\setminus J}
  \overline{\delta}_p(\cF_{J,U}).
\end{equation}
Averaging \eqref{eq:fiber} over all $J\in\binom{[n]}{2d-1}$ and all orderings $\boldsymbol{\pi}$ of $J$, and using \eqref{eq:defect} and \eqref{eq:alpha}, we know that to prove $|\cF|\leq\binom{n}{d}$, it is enough to show that $\alpha_{2d-1}(\cF)\geq0$. Indeed, we could also determine all families attaining equality via our approach.

A family on a specified finite set $J$ is called \emph{elementary} if it is empty, consists of one set, or equals $\binom{J}{s}$ for some $0\le s\le |J|$. We call all other families \emph{exceptional}. For $\cA\subseteq 2^{J}$ with $|J|=p$, put $e_{p}(\cA)=0$ if $\cA$ is elementary and $e_{p}(\cA)=1$ otherwise. In particular, $\cA\subsetneq\binom{J}{s}$ is exceptional whenever $|\cA|\ge 2$.
Corresponding to \eqref{eq:alpha}, for odd $3\le p\le n$, we also define
\begin{equation}\label{eq:beta}
    \beta_p(\cF)
    =\frac{1}{\binom{n}{p}}
  \sum_{J\in\binom{[n]}{p}}
  \sum_{U\subseteq[n]\setminus J}
  e_p(\cF_{J,U}).
\end{equation}

Note the complete layer $\binom{[n]}{d}$ attains the desired bound $\binom{n}{d}$, and every nonempty restricted family of the complete layer $\binom{[n]}{d}$ must be a complete layer. Therefore we will show $\beta_{2d-1}(\cF)=0$ when $\alpha_{2d-1}(\cF)=0$.

For the fixed family $\cF$, write $\alpha_p=\alpha_p(\cF)$ and
$\beta_p=\beta_p(\cF)$. However, it is hard to prove $\alpha_{2d-1}\ge 0$ directly. Therefore, we will firstly prove $\alpha_3\ge0$ and $\beta_3\le3\alpha_3$ in Lemma~\ref{lem:base}. Starting from $\alpha_3\ge0$ and
$\beta_3\le3\alpha_3$, we repeatedly lift $J$ by two coordinates to prove $\alpha_q\ge0$ and
$\beta_q\le M_q\alpha_q$ for every odd $5\le q\le2d-1$ when $n$ is sufficiently large in terms of $d$. Here $M_q$ is a positive number only depending on $q$ when $n$ and $d$ are both fixed. We will show that $M_q$ can be chosen as $(q-3)^{6}2^{4q-16}$ in Proposition~\ref{Proposition:block}.

We firstly prove that for every ordering $\boldsymbol{\pi}$, $\delta_{p,\boldsymbol{\pi}}$ vanishes on every elementary family. On three coordinates, $\overline{\delta}_{3}$ is strictly positive for every exceptional Sperner family.

\begin{lemma}\label{lem:base}
Every elementary family $\cA$ on a ground set $J$ of odd size $p\ge3$ satisfies $\delta_{p,\boldsymbol{\pi}}(\cA)=0$ for every ordering $\boldsymbol{\pi}$ of $J$, and hence $\overline{\delta}_p(\cA)=0$.
Every Sperner family $\cA$ on a three-element ground set satisfies $\overline{\delta}_3(\cA)\ge0$ and $e_3(\cA)\le3\overline{\delta}_3(\cA)$.
\end{lemma}

\begin{proof}[Proof of Lemma~\ref{lem:base}]
Recall that $\delta_{p,\boldsymbol{\pi}}(\cA)=|\Sh(\cA)\cap\cB_{p}(\boldsymbol{\pi})|-|\cA|$ and $\overline{\delta}_{p}(\cA)=\E_{\boldsymbol{\pi}}\delta_{p,\boldsymbol{\pi}}(\cA)$.
If $\cA$ is the empty family, then $\delta_{p,\boldsymbol{\pi}}(\cA)=0$ for any ordering $\boldsymbol{\pi}$ since the empty family can not shatter any set. If $\cA$ is a family consisting of only one set, then $\cA$ can only shatter the empty set. Hence $\overline{\delta}_p(\cA)=0$ for these two cases.

Now we suppose $\cA$ is a complete $s$ layer with $1\le s\le p$. Then $A\subseteq J$ is shattered by $\cA$ if and only if $|A|\le\min(s,p-s)$. In fact, the empty trace needs $\cA$ to have a set disjoint with $A$, and the full trace needs $\cA$ to have a set containing $A$. So $|A|\le s$ and $|A|\le p-s$. Conversely, suppose $|A|\le\min(s,p-s)$. For any subset of $A$, we can add some elements in $J\setminus A$ to this subset to get a set of size $s$. Hence by Lemma~\ref{lem:ballot}, we have $|\Sh(\cA)\cap\cB_p(\boldsymbol{\pi})|=\sum_{k=0}^{\min(s,p-s)}b_{p,k}=\binom{p}{s}=|\cA|$ for every ordering $\boldsymbol{\pi}$. So $\overline{\delta}_p(\cA)=0$ if $\cA$ is a complete layer.

Now we suppose $\cA$ is a Sperner family on a ground set $J$ with $|J|=3$. If $\cA$ is exceptional, then $|\cA|$ must be $2$ since $\cA$ is a Sperner family. If $\cA$ is contained in a complete $1$ layer, then $\cA$ can only shatter the empty set and the sets contained in $\cA$. Note $b_{3,0}=1$ and $b_{3,1}=2$. Then by \eqref{eq:defect}, we have $\overline{\delta}_3(\cA)=\frac{1}{3}$. In the same way, if $\cA$ is contained in a complete $2$ layer, then $\overline{\delta}_3(\cA)=\frac{1}{3}$. If $\cA$ consists of one singleton and one  two-element set, then these two sets must be disjoint since $\cA$ is Sperner. Then $\cA$ can shatter the empty set and all singletons. So $\overline{\delta}_3(\cA)=1$. If $\cA$ is elementary, then $e_3(\cA)=0$ and $\overline{\delta}_3(\cA)=0$. Hence $\overline{\delta}_3(\cA)\ge0$ and $e_3(\cA)\le3\overline{\delta}_3(\cA)$ for any Sperner family $\cA$ on a three-element ground set. 
\end{proof}

We now lift $J$ by two coordinates. Ideally, we would like to get $\alpha_r\le\alpha_{r+2}$, but this does not hold in general. The possible loss comes from the fact that adding a subset of the two extra coordinates to a ballot set still gives a ballot set, unless it already has the largest allowed size and both extra coordinates are added. The following inequality describes this possible loss.

For disjoint finite sets $J,W$ and orderings $\boldsymbol{\pi}$ of $J$ and $\boldsymbol{\rho}$ of $W$, let $\boldsymbol{\pi}\boldsymbol{\rho}$ denote the ordering of $J\cup W$ obtained by listing the elements of $J$ in the order $\boldsymbol{\pi}$, followed by the elements of $W$ in the order $\boldsymbol{\rho}$.

\begin{lemma}\label{lem:ordered}
Let $r=2t+1\ge3$ and $K=J\cup W$, where $J\cap W=\emptyset$, $|J|=r$ and $|W|=2$. Fix orderings $\boldsymbol{\pi}$ of $J$ and $\boldsymbol{\rho}$ of $W$.
For $\cA\subseteq 2^{K}$, put $\cR_{S}(\cA)=\{T\subseteq W:S\in\Sh(\cA_{J,T})\}$ for $S\subseteq J$. Then
\[
\delta_{r+2,\boldsymbol{\pi}\boldsymbol{\rho}}(\cA)\ge\sum_{T\subseteq W}\delta_{r,\boldsymbol{\pi}}(\cA_{J,T})-\bigl|\{S\in\cB_{r}(\boldsymbol{\pi}):|S|=t,\cR_{S}(\cA)=2^{W}\}\bigr|.
\]
\end{lemma}

\begin{proof}[Proof of Lemma~\ref{lem:ordered}]
As in the proof of Lemma~\ref{lem:residual}, if we fix $S\subseteq J$, then for any $T\in\Sh(\cR_S(\cA))$, we have $S\cup T\in\Sh(\cA)$.
Now we fix $S\in\cB_r(\boldsymbol{\pi})$ with $|S|=k$, where $0\le k\le t$.
If $k<t$, then for any $T\subseteq W$, we have $S\cup T\in\cB_{r+2}(\boldsymbol{\pi}\boldsymbol{\rho})$, since $r+1\ge2(k+1)$ and $r+2\ge2(k+2)$.
In the same way, if $k=t$ and $T\subsetneq W$, we also have $S\cup T\in\cB_{r+2}(\boldsymbol{\pi}\boldsymbol{\rho})$.
So we have
\[
|\Sh(\cA)\cap\cB_{r+2}(\boldsymbol{\pi}\boldsymbol{\rho})|\ge\sum_{\substack{S\in\cB_r(\boldsymbol{\pi})\\|S|<t}}|\Sh(\cR_S(\cA))|+\sum_{\substack{S\in\cB_r(\boldsymbol{\pi})\\|S|=t}}|\Sh(\cR_S(\cA))\cap\{T\subseteq W:|T|\le1\}|.
\]
Note that for any $\cR\subseteq2^W$, we have $|\Sh(\cR)|\ge|\cR|$ and $|\Sh(\cR)\cap\{T\subseteq W:|T|\le1\}|$ is at least $|\cR|-\ind_{\{\cR=2^W\}}$.
Here $\ind$ is the characteristic function: $\ind_{\{\cR=2^W\}}=1$ if and only if $\cR=2^W$, and $\ind_{\{\cR=2^W\}}=0$ otherwise.
So we have
\[
|\Sh(\cA)\cap\cB_{r+2}(\boldsymbol{\pi}\boldsymbol{\rho})|\ge\sum_{S\in\cB_r(\boldsymbol{\pi})}|\cR_S(\cA)|-\sum_{\substack{S\in\cB_r(\boldsymbol{\pi})\\|S|=t}}\ind_{\{\cR_S(\cA)=2^W\}}.
\]
By the definition of $\cR_S(\cA)$, we know $|\cR_S(\cA)|=\sum_{T\subseteq W}\ind_{\{S\in\Sh(\cA_{J,T})\}}$, therefore
\begin{align*}
\sum_{S\in\cB_r(\boldsymbol{\pi})}|\cR_S(\cA)|
&=\sum_{S\in\cB_r(\boldsymbol{\pi})}\sum_{T\subseteq W}\ind_{\{S\in\Sh(\cA_{J,T})\}}\\
&=\sum_{T\subseteq W}\sum_{S\subseteq J}\ind_{\{S\in\Sh(\cA_{J,T})\cap\cB_r(\boldsymbol{\pi})\}}\\
&=\sum_{T\subseteq W}|\Sh(\cA_{J,T})\cap\cB_r(\boldsymbol{\pi})|.
\end{align*}
Note that $|\cA|=\sum_{T\subseteq W}|\cA_{J,T}|$ and recall $\delta_{r+2,\boldsymbol{\pi}\boldsymbol{\rho}}(\cA)=|\Sh(\cA)\cap\cB_{r+2}(\boldsymbol{\pi}\boldsymbol{\rho})|-|\cA|$, so we have
\[
\delta_{r+2,\boldsymbol{\pi}\boldsymbol{\rho}}(\cA)\ge\sum_{T\subseteq W}\delta_{r,\boldsymbol{\pi}}(\cA_{J,T})-|\{S\in\cB_r(\boldsymbol{\pi}):|S|=t,\cR_S(\cA)=2^W\}|.
\]

\end{proof}

Now we will average the inequality in Lemma~\ref{lem:ordered} to show the inequality between $\alpha_r$ and $\alpha_{r+2}$. Before this, we need some lemma to handle the error term $\bigl|\{S\in\cB_{r}(\boldsymbol{\pi}):|S|=t,\cR_{S}(\cA)=2^{W}\}\bigr|$.

Let $\cH\subseteq2^{[m]}$ and $i\in[m]$.
For each $H\in\cH$, define
\[
H^{\downarrow i}=
\begin{cases}
H\setminus\{i\},&\text{if }i\in H\text{ and }H\setminus\{i\}\notin\cH,\\
H,&\text{otherwise}.
\end{cases}
\]
The \emph{downshift} of $\cH$ at coordinate $i$ is the family $\cH^{\downarrow i}=\{H^{\downarrow i}:H\in\cH\}$. A \emph{coordinate square} is $Q(B,W)=\{B\cup T:T\subseteq W\}$, where $|W|=2$ and $B\cap W=\emptyset$. We write $c_{2}(\cH)$ for the number of pairs $(B,W)$ such that $Q(B,W)\subseteq\cH$.

The following lemma bounds the number of coordinate squares. 

\begin{lemma}\label{lem:downshift}
Let \(a\ge 1\) be an integer and $\cH\subseteq 2^{[m]}$ be a set system with VC-dimension at most $a$. Then $c_{2}(\cH)\le\binom{a}{2}|\cH|$.
\end{lemma}

\begin{proof}[Proof of Lemma~\ref{lem:downshift}]
 By the definition of downshift, we know $|\cH|=|\cH^{\downarrow i}|$. Now suppose $S\in\Sh(\cH^{\downarrow i})$. If $i\in S$, then for any $T\subseteq S\setminus\{i\}$, we have $H\in \cH^{\downarrow i}$ such that $H\cap S=T\cup \{i\}$. By the definition of downshift, we know $H$ and $H\setminus\{i\}$ are both contained in $\cH$, which means $\cH$ shatters $S$. If $i\notin S$, then for any $T\subseteq S$, we have $H\in\cH^{\downarrow i}$ such that $H\cap S=T$. $\cH$ shatters $S$ since at least one of $H\setminus\{i\}$ and $H\cup\{i\}$ is contained in $\cH$. So $\Sh(\cH^{\downarrow i})\subseteq\Sh(\cH)$.

\begin{claim}\label{clm:shiftsquares}
A coordinate downshift does not decrease the number of coordinate squares.
\end{claim}

\begin{poc}
If a coordinate square $Q(B,W)$ with $W=\{i,j\}$ is contained in $\cH$, then $Q(B,W)$ is contained in $\cH^{\downarrow i}$ after the downshifting of $\cH$ at $i$.
Now we suppose $i\notin B\cup W$. If both $Q(B,W)$ and $Q(B\cup\{i\},W)$ are contained in $\cH$, then any member in $Q(B\cup\{i\},W)$ can not delete $i$ since $Q(B,W)$ is contained in $\cH$. If only one of the $Q(B,W)$ and $Q(B\cup\{i\},W)$ is contained in $\cH$, then we have $Q(B,W) \subseteq\cH^{\downarrow i}$. In fact, if $Q(B\cup\{i\},W)\subseteq\cH$ and $Q(B,W)\nsubseteq\cH$, then any member in $Q(B,W)$ but not in $\cH$ can be obtained after the downshifting of $\cH$ at $i$. The other case is obvious. Note every coordinate square with $i\notin W$ belongs to exactly one such pair, so the number of coordinate squares can not decrease after downshifting.
\end{poc}

Note if $A$ contains $i$ and belongs to $\cH^{\downarrow i}$, then we can not delete $i$ in $A$ since $A\setminus\{i\}$ is already in $\cH^{\downarrow i}$. So after the downshifting of $\cH$ at $i$, we can not delete $i$ in any set containing $i$. Moreover, a downshift at another coordinate preserves this property. 
Let $\cC$ be the family obtained by downshifting $\cH$ at each coordinate once for any ordering.
So we can not delete any coordinate in $\cC$. We also have $|\cC|=|\cH|$ and $\Sh(\cC)\subseteq\Sh(\cH)$. Any set in $\cC$ can be shattered by $\cC$ since every subset of this set is already contained in $\cC$. So any set in $\cC$ has size at most $a$. Note a coordinate square $Q(B,W)$ is contained in $\cC$ if and only if $B\cup W$ is contained in $\cC$. Hence by  Claim~\ref{clm:shiftsquares}, we have $c_{2}(\cH)\le c_{2}(\cC)=\sum_{A\in\cC}\binom{|A|}{2}\le\binom{a}{2}|\cC|=\binom{a}{2}|\cH|$.

\end{proof}

Suppose that $|J|=2t+1$ and $|S|=t$. We will show that if $U$ belongs to a coordinate square in $\cR_{J,S}$, then $\cF_{J,U}$ is exceptional, so that the number of coordinate squares can be bounded by $\beta$.

\begin{lemma}\label{lemma:squares}
Let $J\subseteq[n]$ have size $2t+1\ge3$, where $t\le d$, and let $S\in\binom{J}{t}$.
Put $\cP_{J,S}=\{U\in\cR_{J,S}:e_{2t+1}(\cF_{J,U})=1\}$.
Every coordinate square contained in $\cR_{J,S}$ is also contained in $\cP_{J,S}$.
Consequently,
\[
c_{2}(\cR_{J,S})=c_{2}(\cP_{J,S})\le\binom{d-t}{2}|\cP_{J,S}|.
\]
Moreover, $\sum_{S\in\binom{J}{t}}|\cP_{J,S}|\le\binom{2t+1}{t}\sum_{U\subseteq[n]\setminus J}e_{2t+1}(\cF_{J,U})$.
\end{lemma}

\begin{proof}[Proof of Lemma~\ref{lemma:squares}]
By the proof of Lemma~\ref{lem:base}, we know if $\cF_{J,U}$ is elementary and shatters $S$, then $\cF_{J,U}$ must be $\binom{J}{t}$ or $\binom{J}{t+1}$.

\begin{claim}\label{clm:residualchain}
If $U_0\subsetneq U_1\subsetneq U_2$ belong to $\cR_{J,S}$, then all three restricted families $\cF_{J,U_i}$ are exceptional.
\end{claim}

\begin{poc}
Note if $U\subsetneq V$, then any set in $\cF_{J,U}$ can not be contained in any set in $\cF_{J,V}$ since $\cF$ is a Sperner family. We write $\cA_i=\cF_{J,U_i}$ for $i=0,1,2$.

Now suppose $\cA_1=\binom{J}{t}$. Then any set in $\cA_2$ has size at most $t-1$. Hence $\cA_2$ can not shatter $S$, which is contradictory. In the same way, $\cA_0$ can not be $\binom{J}{t}$. If $\cA_1=\binom{J}{t+1}$, then any set in $\cA_0$ has size at least $t+2$. But the empty trace needs there to exist a set in $\cA_0$ such that this set is disjoint with $S$. Note $|J|=2t+1$ and $|S|=t$. This is impossible. In the same way, $\cA_2$ can not be $\binom{J}{t+1}$. 

Now if $\cA_0=\binom{J}{t+1}$ or $\cA_2=\binom{J}{t}$, note $\cA_0=\binom{J}{t+1}$ means $S$ is contained in $\cA_1$ since $\cA_1$ shatters $S$ and any set in $\cA_1$ has size at most $t$. Note $\cA_2$ also shatters $S$, which means $S$ is contained in some set in $\cA_2$. This is contradictory. Similarly, $\cA_2=\binom{J}{t}$ means $J\setminus S$ is contained in $\cA_1$, which is contradictory since $\cA_0$ shatters $S$. Therefore, all three families are exceptional.
\end{poc}

Take any coordinate square $Q(B,W)$ in $\cR_{J,S}$. Then every member of $Q(B,W)$ is contained in a chain $U_0\subsetneq U_1 \subsetneq U_2$. By Claim~\ref{clm:residualchain} and the definition of $\cP_{J,S}$, we know $Q(B,W)$ is contained in $\cP_{J,S}$. Hence $c_{2}(\cP_{J,S})=c_{2}(\cR_{J,S})$ and by Lemma~\ref{lem:residual}, we have $\VC(\cP_{J,S})\le d-t$. If $d>t$, then by Lemma~\ref{lem:downshift}, we have $c_{2}(\cP_{J,S})\le\binom{d-t}{2}|\cP_{J,S}|$. If $d=t$ and $\cP_{J,S}$ contains a coordinate square $Q(B,W)$, then $\cP_{J,S}$ can shatter the two-element set $W$, which contradicts $\VC(\cP_{J,S})\le d-t$. So if $d=t$, we must have $c_{2}(\cR_{J,S})=c_{2}(\cP_{J,S})=0$. Finally, for a fixed $S$, $|\cP_{J,S}|$ is the number of restricted families $\cF_{J,U}$ such that $\cF_{J,U}$ shatters $S$ and $e_{2t+1}(\cF_{J,U})=1$. Note $\cF_{J,U}$ can shatter at most $\binom{2t+1}{t}$ sets of size $t$. Then we have the last inequality.
\end{proof}

For odd $r$ with $3\le r\le\min(2d+1,n-2)$, put $a_{r}(n,d)=b_{r,\frac{(r-1)}{2}}\cdot\frac{\binom{d-(r-1)/2}{2}}{\binom{n-r}{2}}$. Now we average the preceding inequality in Lemma~\ref{lem:ordered} and show that $\alpha_{r+2}\ge \alpha_{r}-a_{r}(n,d)\beta_{r}$.

\begin{prop}\label{prop:transport}
For every odd $r$ with $3\le r\le\min(2d+1,n-2)$, we have
\[
\alpha_{r+2}\ge \alpha_{r}-a_{r}(n,d)\beta_{r}.
\]
\end{prop}

\begin{proof}[Proof of Proposition~\ref{prop:transport}]
Put $t=\frac{r-1}{2}$.
We first choose $J\in\binom{[n]}{r}$ uniformly, and then choose $W\in\binom{[n]\setminus J}{2}$ uniformly.
Further, we choose an ordering $\boldsymbol{\pi}$ of $J$ and an ordering $\boldsymbol{\rho}$ of $W$ uniformly and independently.
Let $K=J\cup W$.
This is equivalent to choosing $K\in\binom{[n]}{r+2}$ uniformly and choosing an ordering $\boldsymbol{\sigma}$ of $K$ uniformly.
Now we fix such $J,W,\boldsymbol{\pi},\boldsymbol{\rho}$, and also fix $Z\subseteq[n]\setminus K$.
By Lemma~\ref{lem:ordered}, we have
\[
\delta_{r+2,\boldsymbol{\pi}\boldsymbol{\rho}}(\cF_{K,Z})\ge\sum_{T\subseteq W}\delta_{r,\boldsymbol{\pi}}((\cF_{K,Z})_{J,T})-|\{S\in\cB_r(\boldsymbol{\pi}):|S|=t,\cR_S(\cF_{K,Z})=2^W\}|.
\]
Summing over $Z\subseteq[n]\setminus K$ and averaging over $J,W,\boldsymbol{\pi},\boldsymbol{\rho}$, we have $f_1\ge f_2-f_3$, where
\begin{align*}
f_1&=\frac{1}{\binom{n}{r}\binom{n-r}{2}r!2!}\sum_{J\in\binom{[n]}{r}}\sum_{W\in\binom{[n]\setminus J}{2}}\sum_{\boldsymbol{\pi}}\sum_{\boldsymbol{\rho}}\sum_{Z\subseteq[n]\setminus K}\delta_{r+2,\boldsymbol{\pi}\boldsymbol{\rho}}(\cF_{K,Z}),\\
f_2&=\frac{1}{\binom{n}{r}\binom{n-r}{2}r!2!}\sum_{J\in\binom{[n]}{r}}\sum_{W\in\binom{[n]\setminus J}{2}}\sum_{\boldsymbol{\pi}}\sum_{\boldsymbol{\rho}}\sum_{Z\subseteq[n]\setminus K}\sum_{T\subseteq W}\delta_{r,\boldsymbol{\pi}}((\cF_{K,Z})_{J,T}),\\
f_3&=\frac{1}{\binom{n}{r}\binom{n-r}{2}r!2!}\sum_{J\in\binom{[n]}{r}}\sum_{W\in\binom{[n]\setminus J}{2}}\sum_{\boldsymbol{\pi}}\sum_{\boldsymbol{\rho}}\sum_{Z\subseteq[n]\setminus K}|\{S\in\cB_r(\boldsymbol{\pi}):|S|=t,\cR_S(\cF_{K,Z})=2^W\}|.
\end{align*}
Then
\begin{align*}
f_1&=\frac{1}{\binom{n}{r+2}(r+2)!}\sum_{K\in\binom{[n]}{r+2}}\sum_{\boldsymbol{\sigma}}\sum_{Z\subseteq[n]\setminus K}\delta_{r+2,\boldsymbol{\sigma}}(\cF_{K,Z})\\
&=\frac{1}{\binom{n}{r+2}}\sum_{K\in\binom{[n]}{r+2}}\sum_{Z\subseteq[n]\setminus K}\overline{\delta}_{r+2}(\cF_{K,Z})
=\alpha_{r+2}.
\end{align*}
Note that $(\cF_{K,Z})_{J,T}=\{A\subseteq J:A\cup T\cup Z\in\cF\}=\cF_{J,Z\cup T}$, so we have
\[
f_2=\frac{1}{\binom{n}{r+2}(r+2)!}\sum_{J\in\binom{[n]}{r}}\sum_{W\in\binom{[n]\setminus J}{2}}\sum_{\boldsymbol{\pi}}\sum_{\boldsymbol{\rho}}\sum_{Z\subseteq[n]\setminus K}\sum_{T\subseteq W}\delta_{r,\boldsymbol{\pi}}(\cF_{J,Z\cup T}).
\]
Since every $U\subseteq[n]\setminus J$ has a unique decomposition $U=Z\cup T$, with $Z=U\setminus K$ and $T=U\cap(K\setminus J)$, we have
\[
f_2=\frac{1}{\binom{n}{r+2}(r+2)!}\sum_{J\in\binom{[n]}{r}}\sum_{W\in\binom{[n]\setminus J}{2}}\sum_{\boldsymbol{\pi}}\sum_{\boldsymbol{\rho}}\sum_{U\subseteq[n]\setminus J}\delta_{r,\boldsymbol{\pi}}(\cF_{J,U}).
\]
Note that $\delta_{r,\boldsymbol{\pi}}(\cF_{J,U})$ is independent of the choice of $W$ and $\boldsymbol{\rho}$, so we have
\begin{align*}
f_2&=\frac{\binom{n-r}{2}2!r!}{\binom{n}{r+2}(r+2)!}\sum_{J\in\binom{[n]}{r}}\sum_{U\subseteq[n]\setminus J}(\sum_{\boldsymbol{\pi}}\frac{1}{r!}\delta_{r,\boldsymbol{\pi}}(\cF_{J,U}))\\
&=\frac{1}{\binom{n}{r}}\sum_{J\in\binom{[n]}{r}}\sum_{U\subseteq[n]\setminus J}\overline{\delta}_r(\cF_{J,U})
=\alpha_r.
\end{align*}
Now we compute $f_3$.
We firstly prove that $\cR_S(\cF_{K,Z})=2^W$ if and only if $Q(Z,W)\subseteq\cR_{J,S}$.
Recall the definition of $\cR_S(\cA)$ in Lemma~\ref{lem:ordered}: $\cR_{S}(\cA)=\{T\subseteq W:S\in\Sh(\cA_{J,T})\}$ for $S\subseteq J$. Note we have $(\cF_{K,Z})_{J,T}=\cF_{J,T\cup Z}$, so $\cR_S(\cF_{K,Z})=2^W$ if and only if $S\in\Sh((\cF_{K,Z})_{J,T})=\Sh(\cF_{J,T\cup Z})$ for every $T\subseteq W$, which means $T\cup Z\in\cR_{J,S}$. 
Recall that a coordinate square $Q(Z,W)$ is defined as $Q(Z,W)=\{Z\cup T:T\subseteq W\}$.
Hence $\cR_S(\cF_{K,Z})=2^W$ if and only if $Q(Z,W)\subseteq\cR_{J,S}$.
So we have
\[
f_3=\frac{1}{\binom{n}{r}\binom{n-r}{2}r!2!}\sum_{J\in\binom{[n]}{r}}\sum_{W\in\binom{[n]\setminus J}{2}}\sum_{\boldsymbol{\pi}}\sum_{\boldsymbol{\rho}}\sum_{Z\subseteq[n]\setminus K}\sum_{\substack{S\in\cB_r(\boldsymbol{\pi})\\|S|=t}}\ind_{\{Q(Z,W)\subseteq\cR_{J,S}\}}.
\]
Now we fix $J\in\binom{[n]}{r}$ and $S\in\binom{J}{t}$.
For every $W\in\binom{[n]\setminus J}{2}$ and $Z\subseteq[n]\setminus(J\cup W)$, we have a coordinate square $Q(Z,W)$.
So $\sum_{W\in\binom{[n]\setminus J}{2}}\sum_{Z\subseteq[n]\setminus(J\cup W)}\ind_{\{Q(Z,W)\subseteq\cR_{J,S}\}}=c_2(\cR_{J,S})$.
Note that we choose $\boldsymbol{\pi}$ uniformly, so by Lemma~\ref{lem:ballot}, we have $\mathbb{P}(S\in\cB_r(\boldsymbol{\pi}))=b_{r,t}/\binom{r}{t}$. Also note that $c_2(\cR_{J,S})$ is independent of the choice of $\boldsymbol{\rho}$, so we have
\[
f_3=\frac{b_{r,t}}{\binom{n}{r}\binom{n-r}{2}\binom{r}{t}}\sum_{J\in\binom{[n]}{r}}\sum_{S\in\binom{J}{t}}c_2(\cR_{J,S}).
\]
By Lemma~\ref{lemma:squares}, we have
\[
f_3\le\frac{b_{r,t}\binom{d-t}{2}\binom{r}{t}}{\binom{n}{r}\binom{n-r}{2}\binom{r}{t}}\sum_{J\in\binom{[n]}{r}}\sum_{U\subseteq[n]\setminus J}e_r(\cF_{J,U}).
\]
So we have $f_3\le a_r(n,d)\beta_r$, hence
\begin{equation}\label{eq:transport}
\alpha_{r+2}\ge\alpha_r-a_r(n,d)\beta_r.
\end{equation}
\end{proof}

Even if we know $\alpha_r\ge0$, we cannot get $\alpha_{r+2}\ge0$ only by Proposition~\ref{prop:transport}. Therefore we need to dig out more information between $\beta$ and $\alpha$. Before this, we need several auxiliary results.

We first show the whole family $\cF$ is strongly influenced by the restricted family $\cF_{J,U}$.

\begin{lemma}\label{lem:rigidity}
Let $3\le p\le m$, and let $\cA\subseteq2^{[m]}$ be a Sperner family.
Suppose that $\cA_{J,U}$ is elementary for every $J\in\binom{[m]}{p}$ and every $U\subseteq[m]\setminus J$.
Then either $\cA$ is elementary, or $|X\triangle Y|>p$ for all distinct $X,Y\in\cA$.
\end{lemma}

\begin{proof}[Proof of Lemma~\ref{lem:rigidity}]
Induction on $m$. If $m=p$ and there exist two distinct sets $X,Y\in\cA$ such that $|X\triangle Y|\le p$, then $\cA$ must be a complete layer. In fact, in this case, $J$ can only be $[m]$, which means $\cA$ is elementary. Since $\cA$ has at least two sets, $\cA$ must be a complete $s$ layer with $1\le s\le m-1$.

Assume that the statement holds for $m$ with $m\ge p$. We prove it for $m+1$. If there exist two distinct sets $X,Y\in\cA$ such that $|X\triangle Y|\le p$, then there exists at least one coordinate not contained in $X\triangle Y$. So we may suppose $i$ is not contained in either $X$ or $Y$ for some $i$. Indeed, if we can not take such $i$, then there exists $j$, such that $j\in X\cap Y$. So we can take complements for every set in $\cA$. Note taking complements preserves the hypothesis, since the complement in $J$ of an elementary restricted family is again elementary. 

Without loss of generality, we may further assume that the coordinate $i$ we take is $m+1$. So $X,Y\in 2^{[m]}$ and $|X\triangle Y|\le p$. By the induction hypothesis, we know that the subfamily of $\cA$ consisting of all sets not containing $m+1$ is a complete $s$ layer on the ground set $[m]$ with $1\le s\le m-1$. For any $B\in\binom{[m]}{s-1}$, $B\cup \{m+1\}$ is also contained in $\cA$. In fact, since $1\le s\le m-1$, we can take two distinct coordinates $a,b$ outside $B\cup \{m+1\}$. Note $p\ge 3$, so we can choose a set $J$ of size $p$ such that $J$ contains $\{a,b,m+1\}$. We take the outside part to be $B\setminus J$. Then this restricted family contains at least two members and hence it must be a complete layer, which means $B\cup \{m+1\}$ must be contained in $\cA$. So $\cA$ contains a complete $s$ layer on the ground set $[m+1]$. Note $\cA$ is a Sperner family, so no other set can be contained in $\cA$. Hence $\cA$ is a complete $s$ layer on the ground set $[m+1]$ with $1\le s\le m-1$.

\end{proof}

We next compare $\overline{\delta}_{q}(\cA)$ with the number of exceptional restricted families $\cA_{J,U}$ for $|J|=p$.

\begin{lemma}\label{lem:localrepair}
Let $p,q$ be odd with $3\le p<q\le2p-1$.
Every Sperner family $\cA\subseteq2^{[q]}$ satisfies
\[
e_q(\cA)\le\overline{\delta}_q(\cA)+\binom{q}{(q-1)/2}\sum_{J\in\binom{[q]}{p}}\sum_{U\subseteq[q]\setminus J}e_p(\cA_{J,U}).
\]
\end{lemma}

\begin{proof}[Proof of Lemma~\ref{lem:localrepair}]
If there exists a restricted family $\cA_{J,U}$ with $|J|=p$ that is exceptional, then the sum in the statement is at least one. Note $\cA$ is nonempty and hence shatters the empty set. So by \eqref{eq:defect}, we have $\overline{\delta}_q(\cA)\ge1-|\cA|$. Using Lemma~\ref{lem:Sperner}, we have $|\cA|\le\binom{q}{(q-1)/2}$, so the right side of the stated inequality is at least $1-|\cA|+\binom{q}{(q-1)/2}\ge1\ge e_q(\cA)$.

Now we suppose every restricted family $\cA_{J,U}$ with $|J|=p$ is elementary. If $\cA$ is elementary, then $e_q(\cA)=\overline{\delta}_q(\cA)=0$ by Lemma~\ref{lem:base}.
If $\cA$ is exceptional, then $|\cA|\ge2$ and Lemma~\ref{lem:rigidity} gives $|A\triangle B|\ge p+1$ for all distinct $A,B\in\cA$. 

\begin{claim}\label{clm:separatedsize}
 $|\cA|\le\frac{2}{3}(p+1)$.
\end{claim}

\begin{poc}
Choose $A,B$ independently and uniformly from $\cA$, allowing $A=B$.
For a fixed coordinate, let $\theta$ be the proportion of members of $\cA$ containing it.
By independence, $A$ and $B$ differ at this coordinate with probability $2\theta(1-\theta)\le\frac{1}{2}$.
Summing over the $q$ coordinates gives $\E|A\triangle B|\le \frac{q}{2}$.
On the other hand, $\mathbb{P}(A=B)=\frac{1}{|\cA|}$, and $|A\triangle B|\ge p+1$ whenever $A\neq B$.
Thus $(p+1)(1-\frac{1}{|\cA|})\le\E|A\triangle B|\le \frac{q}{2}$.
Since $2(p+1)-q\ge3$, we have $|\cA|\le\frac{2(p+1)}{2(p+1)-q}\le\frac{2}{3}(p+1)$.
\end{poc}

Fix distinct $A,B\in\cA$.
For each $i\in A\triangle B$, one of $A$ and $B$ contains $i$ and the other does not, so $\{i\}$ is shattered by $\cA$.
Thus $\cA$ shatters at least $p+1$ singletons and the empty set.
Note $b_{q,0}=1$ and $b_{q,1}=q-1$, so we have $\overline{\delta}_q(\cA)\ge1+\frac{2}{3}(p+1)-|\cA|\ge1=e_q(\cA)$.

\end{proof}

For these values of $p,q$, put $c_{q,p}=\binom{q}{(q-1)/2}\binom{q}{p}$. Summing Lemma~\ref{lem:localrepair} over the restricted families gives the following inequality.

\begin{prop}\label{prop:repair}
For odd $p,q$ with $3\le p<q\le\min(2p-1,n)$, we have
\[
\beta_{q}\le \alpha_{q}+c_{q,p}\beta_{p}.
\]
\end{prop}

\begin{proof}[Proof of Proposition~\ref{prop:repair}]
Fix $K\in\binom{[n]}{q}$ and $Z\subseteq[n]\setminus K$.
For $J\in\binom{K}{p}$ and $T\subseteq K\setminus J$, we have $(\cF_{K,Z})_{J,T}=\{A\subseteq J:A\cup T\cup Z\in\cF\}=\cF_{J,Z\cup T}$.
For fixed $J\subseteq K$, every $U\subseteq[n]\setminus J$ has a unique decomposition $U=Z\cup T$, with $Z=U\setminus K$ and $T=U\cap(K\setminus J)$.
Hence we have
\[
\sum_{Z\subseteq[n]\setminus K}\sum_{T\subseteq K\setminus J}e_p(\cF_{J,Z\cup T})=\sum_{U\subseteq[n]\setminus J}e_p(\cF_{J,U}).
\]
Applying Lemma~\ref{lem:localrepair} to $\cF_{K,Z}$, we have
\[
e_q(\cF_{K,Z})\le\overline{\delta}_q(\cF_{K,Z})+\binom{q}{(q-1)/2}\sum_{J\in\binom{K}{p}}\sum_{T\subseteq K\setminus J}e_p((\cF_{K,Z})_{J,T}).
\]
Summing over $Z\subseteq[n]\setminus K$ and averaging over $K\in\binom{[n]}{q}$, we have
\begin{align*}
\beta_q
&\le\alpha_q+\frac{\binom{q}{(q-1)/2}\binom{n}{p}}{\binom{n}{q}}\frac{1}{\binom{n}{p}}\sum_{K\in\binom{[n]}{q}}\sum_{Z\subseteq[n]\setminus K}\sum_{J\in\binom{K}{p}}\sum_{T\subseteq K\setminus J}e_p((\cF_{K,Z})_{J,T})\\
&=\alpha_q+\frac{\binom{q}{(q-1)/2}\binom{n}{p}}{\binom{n}{q}}\frac{1}{\binom{n}{p}}\sum_{K\in\binom{[n]}{q}}\sum_{J\in\binom{K}{p}}\sum_{U\subseteq[n]\setminus J}e_p(\cF_{J,U}).
\end{align*}
Fix $J$. Then there are $\binom{n-p}{q-p}$ sets $K\in\binom{[n]}{q}$ that contain $J$. Note $\binom{n}{p}\binom{n-p}{q-p}=\binom{n}{q}\binom{q}{p}$, so we have
\begin{equation}\label{eq:repair}
\beta_q\le\alpha_q+\frac{\binom{q}{(q-1)/2}\binom{n}{p}\binom{n-p}{q-p}}{\binom{n}{q}}\beta_p=\alpha_q+c_{q,p}\beta_p.
\end{equation}
\end{proof}

We now suppose that $n$ is sufficiently large in terms of $d$. Combining \eqref{eq:transport} and \eqref{eq:repair}, we will prove that $\alpha_{q}\ge 0$ and $\beta_{q}\le(q-3)^{6}\cdot2^{4q-16}\alpha_{q}$ for every odd $5\le q\le 2d-1$.

\begin{prop}\label{Proposition:block}
For every integer $d\ge3$, put $n_0(d)=d^{10}\cdot2^{5d-10}$.
Then for every $n\ge n_0(d)$ and every odd $q$ with
$5\le q\le2d-1$, we have
\(
\alpha_{q}\ge 0\) and \(\beta_{q}\le(q-3)^{6}\cdot2^{4q-16}\alpha_{q}.
\)
\end{prop}

\begin{proof}[Proof of Proposition~\ref{Proposition:block}]
For any fixed odd $q$ with $5\le q\le2d-1$, put $p=2\lfloor\frac{q+1}{4}\rfloor+1$. Then $p$ is the smallest odd integer at least $\frac{q+1}{2}$. 
For fixed $d$ and any odd $r$ with $p\le r<q$, $a_{r}(n,d)$ tends to zero as $n$ tends to infinity. For this fixed $d$, if $n\ge d^{10}\cdot2^{5d-10}$, then for any such $p,q,r$, we have 
\begin{equation}\label{eq:smallloss}
(1+2^{2q-6}\cdot p^{6}\cdot2^{4p-16})\sum_{\substack{p\le r<q\\ r \text{ odd}}}a_{r}(n,d)\le\frac{1}{2}.
\end{equation}
In fact, recall that $a_{r}(n,d)=\frac{b_{r,\frac{r-1}{2}}\cdot\binom{d-(r-1)/2}{2}}{\binom{n-r}{2}}$ and $b_{r,\frac{r-1}{2}}=\binom{r}{(r-1)/2}-\binom{r}{(r-3)/2}$. Since $5\le q\le2d-1$ and $p=2\lfloor\frac{q+1}{4}\rfloor+1$, we have $3\le p\le\frac{q+3}{2}$ and $\frac{q+3}{2}\le d+1$. Hence we have $r\le q-2\le2d-3$ and $1\le\frac{r-1}{2}\le d-2$.

Since $2\binom{r}{(r-1)/2}\le2^r$, we have $0\le b_{r,\frac{r-1}{2}}\le\binom{r}{(r-1)/2}\le2^{r-1}\le2^{2d-4}$. We also have $\binom{n-r}{2}\ge\frac{n^2}{8}$ since $n\ge d^{10}\cdot2^{5d-10}\ge4d$ and $r\le2d-3$. Note $\binom{d-(r-1)/2}{2}\le\frac{d^2}{2}$, so for any odd $r$ with $p\le r<q$, we have 
\begin{equation}\label{eq:8}
    a_{r}(n,d)\le\frac{2^{2d-4}\cdot d^2/2}{n^2/8}=\frac{d^2\cdot2^{2d-2}}{n^2}.
\end{equation}
Since $q\ge5$ and $p\ge3$, we have $2^{2q-6}\cdot p^{6}\cdot2^{4p-16}=p^6\cdot2^{2q+4p-22}\ge1$. Note $q\le2d-1$ and $p\le d+1$, so we have
\begin{equation}\label{eq:9}
    1+2^{2q-6}\cdot p^{6}\cdot2^{4p-16}\le2^{2q-5}\cdot p^{6}\cdot2^{4p-16}\le (d+1)^6\cdot2^{8d-19}.
\end{equation}
Combining \eqref{eq:8} and \eqref{eq:9}, we have $(1+2^{2q-6}\cdot p^{6}\cdot2^{4p-16})\sum\limits_{\substack{p\le r<q\\ r \text{ odd}}}a_{r}(n,d)\le\frac{d^3(d+1)^6\cdot2^{10d-21}}{n^2}$. So we get \eqref{eq:smallloss} when $n\ge d^{10}\cdot2^{5d-10}$.
In particular, $0\le a_{r}(n,d)\le 1$ for every odd $r$ with $p\le r<q$ when $n\ge d^{10}\cdot2^{5d-10}$.

Now we suppose $d\ge3$ and $n\ge d^{10}\cdot2^{5d-10}$ are both fixed. Suppose $\alpha_{p}\ge 0$ and $\beta_{p}\le (p-3)^{6}\cdot2^{4p-16}\alpha_{p}$ if $p\ge5$. Note $\beta_{3}\le3\alpha_{3}$ by Lemma~\ref{lem:base}, so in any case, we always have $\beta_p\le p^6 \cdot2^{4p-16}\alpha_p$.

\begin{claim}\label{clm:blockloss}
 $\alpha_{q}\ge\frac{1}{2}\alpha_{p}$.
\end{claim}

\begin{poc}
Write $a_{r}=a_{r}(n,d)$. Recall that $c_{r,p}=\binom{r}{(r-1)/2}\binom{r}{p}$, so for any odd $r$ with $p<r<q$, we have $c_{r,p}\le 2^{2r-2}\le 2^{2q-6}$ and hence have $\beta_{r}\le\alpha_{r}+2^{2q-6}\beta_{p}$ by \eqref{eq:repair}. We also have $\beta_{p}\le\alpha_{p}+2^{2q-6}\beta_{p}$ since $\alpha_{p},\beta_{p}\ge 0$ and $2^{2q-6}\ge 1$. Substituting this inequality into 
\eqref{eq:transport}, we have $\alpha_{r+2}+2^{2q-6}\beta_{p}\ge(1-a_{r})(\alpha_{r}+2^{2q-6}\beta_{p})$.
Since $0\le a_r\le1$ and $\alpha_p+2^{2q-6}\beta_p\ge0$, we have
\[
\alpha_q+2^{2q-6}\beta_p\ge(\prod_{\substack{p\le r<q\\ r\text{ odd}}}(1-a_r))(\alpha_p+2^{2q-6}\beta_p).
\]
Using $\prod\limits_{\substack{p\le r<q\\ r\text{ odd}}}(1-a_r)\ge1-\sum\limits_{\substack{p\le r<q\\ r\text{ odd}}}a_r$, we have $\alpha_q
\ge\alpha_p-(\alpha_p+2^{2q-6}\beta_p)\sum\limits_{\substack{p\le r<q\\ r\text{ odd}}}a_r$.
By the hypothesis $\beta_p\le p^6 \cdot2^{4p-16}\alpha_p$, we have 
\[
\alpha_q\ge\left(1-(1+2^{2q-6}\cdot p^{6}\cdot2^{4p-16})\sum_{\substack{p\le r<q\\ r \text{ odd}}}a_{r}(n,d)\right)\alpha_p.
\]
By \eqref{eq:smallloss}, we have $\alpha_q\ge\frac12\alpha_p$.
\end{poc}

Since $\alpha_p\ge0$ and $\alpha_q\ge\frac12\alpha_p$, we have $\alpha_q\ge0$. On the other hand, using \eqref{eq:repair}, we have $\beta_q\le\alpha_q+c_{q,p}\beta_p$.
If $q=5$, then $p=3$ and $c_{5,3}=100$.
Note we have $\alpha_3\ge0$ and $\beta_3\le3\alpha_3$ by Lemma~\ref{lem:base}.
Using $\beta_3\le3\alpha_3$ and $\alpha_5\ge\frac{1}{2}\alpha_3$, we have $\beta_5\le\alpha_5+100\beta_3\le601\alpha_5\le2^{10}\alpha_5$, as required.
If $q\ge7$, then $p\ge5$.
Since $\beta_p\le(p-3)^6 \cdot2^{4p-16}\alpha_p$ and $\alpha_{p}\le2\alpha_{q}$, we have $\beta_q\le(1+2c_{q,p}\cdot(p-3)^6 \cdot2^{4p-16})\alpha_q$.
Note that $q-3\ge2(p-3)$ and $p\le\frac{q+3}{2}$, so
$\frac{(q-3)^6 \cdot2^{4q-16}}{(p-3)^6\cdot 2^{4p-16}}\ge2^{6+4(q-p)}\ge2^{2q}$.
Also note that $c_{q,p}\le2^{2q-2}$ and $c_{q,p}\cdot(p-3)^6 \cdot2^{4p-16}\ge1$, so we have 
\begin{align*}
1+2c_{q,p}\cdot(p-3)^6\cdot 2^{4p-16}
&\le3c_{q,p}\cdot(p-3)^6\cdot 2^{4p-16}\\
&\le2^{2q}\cdot(p-3)^6\cdot 2^{4p-16}\\
&\le(q-3)^6\cdot 2^{4q-16},
\end{align*}
which means $\beta_q\le(q-3)^6\cdot 2^{4q-16}\alpha_q$.

For every odd $q\ge5$, $p=2\lfloor\frac{q+1}{4}\rfloor+1$ is the smallest odd integer at least $\frac{q+1}{2}$, and hence $2\lfloor\frac{p+1}{4}\rfloor+1$ is the smallest odd integer at least $\frac{p+1}{2}$. Repeating this finitely many times, we always get $3$. By Lemma~\ref{lem:base}, we have $\alpha_3\ge0$ and $\beta_3\le3\alpha_3$. So, by induction, if $d$ is fixed and $n\ge d^{10}\cdot2^{5d-10}$, then for every odd $q$ with $5\le q\le 2d-1$, we have $\alpha_{q}\ge 0$ and $\beta_{q}\le(q-3)^{6}\cdot 2^{4q-16}\alpha_{q}$.

\end{proof}

It is enough to prove that $\alpha_{2d-1}(\cF)\geq0$ to get $|\cF|\leq\binom{n}{d}$. So we have proved the bound in the main theorem~\ref{thm:main}. The rest of the equality case can be handled by Lemma~\ref{lem:rigidity}.

\begin{proof}[Proof of Theorem~\ref{thm:main}]
 By Proposition~\ref{Proposition:block}, \eqref{eq:fiber} and \eqref{eq:alpha}, we have $0\le\alpha_{2d-1}\le\binom{n}{d}-|\cF|$ when $n\ge d^{10}\cdot2^{5d-10}$. This proves the bound.

Now we suppose that $|\cF|=\binom{n}{d}$.
Then by Proposition~\ref{Proposition:block}, we know that $\alpha_{2d-1}=0$ means that $\beta_{2d-1}=0$.
Note that $\beta_{2d-1}$ is the average number of exceptional restricted families, so $\beta_{2d-1}=0$ means that every restricted family is elementary.
In other words, for any $J\in\binom{[n]}{2d-1}$ and any $U\subseteq[n]\setminus J$, $\cF_{J,U}$ is elementary.

If $\cF$ is exceptional, then $\cF$ has at least two members.
By Lemma~\ref{lem:rigidity}, for any two distinct sets $X,Y\in\cF$, we have $|X\triangle Y|>2d-1$.
Now we fix $J\in\binom{[n]}{2d-1}$ and denote $\cH=\{A\setminus J:A\in\cF\}\subseteq2^{[n]\setminus J}$.
Note that if $A\setminus J=B\setminus J$ for some $A,B\in\cF$, then we must have $A=B$.
In fact, if $A\neq B$ and $A\setminus J=B\setminus J$, then $|A\triangle B|>2d-1$ and $A\triangle B\subseteq J$, but $|J|=2d-1$, which is impossible.
So we have a bijection between $\cF$ and $\cH$.
Note that $\cH$ is a family on the ground set $[n]\setminus J$, so any set shattered by $\cH$ can also be shattered by $\cF$.
Hence $\VC(\cH)\le\VC(\cF)\le d$.
Then by Lemma~\ref{lem:sauer}, we have $|\cH|\le\sum_{j=0}^{d}\binom{n-2d+1}{j}$.
Note that $b_{2d-1,0}=1$ and $b_{2d-1,1}=2d-2>0$, so by Lemma~\ref{lem:terminal}, we have $\binom{n}{d}=\sum_{s=0}^{d-1}b_{2d-1,s}\sum_{j=0}^{d-s}\binom{n-2d+1}{j}>\sum_{j=0}^{d}\binom{n-2d+1}{j}$, which contradicts $\binom{n}{d}=|\cF|=|\cH|$.

Hence $\cF$ must be elementary. Note that $|\cF|=\binom{n}{d}>1$, so it must be a complete $s$ layer with $1\le s \le n-1$. Then obviously, $s$ must be $d$ or $n-d$. This finishes the proof.
\end{proof}

\section*{Acknowledgement}
The authors would like to thank Jian Wang, Yang Shu and Xiaochen Zhao for many helpful discussions.

\bibliographystyle{abbrv}
\bibliography{bib}

\end{document}